\documentclass[12pt,reqno]{amsart}

\usepackage{epsfig}
\usepackage{color}
\usepackage{amsmath, amsthm, amsfonts, amssymb, mathrsfs,mathtools,faktor}
\usepackage{graphicx}
\usepackage{enumerate}

\numberwithin{equation}{section}

\newcommand{\C}{{\mathbb C}}

\newtheorem{theo}{{\sc \bf Theorem}}[section]

\newtheorem{prop}[theo]{{\sc \bf Proposition}}

\theoremstyle{definition}
\newtheorem{defin}{Definition}[section]

\begin{document}
	
	\title{Compact Linear Relations and Their Spectral Properties}

	\author[Acharya]{Keshav Raj Acharya}
	\address{Department of Mathematics,
		Embry-Riddle Aeronautical University,
		600 S. Clyde Moris Blvd., Daytona Beach, FL 32114, U.S.A.}
	\email{acharyak@erau.edu}
	
	\author[McBride]{Matt McBride}
	\address{Department of Mathematics and Statistics,
		Mississippi State University,
		175 President's Cir., Mississippi State, MS 39762, U.S.A.}
	\email{mmcbride@math.msstate.edu}

	\date{\today}
	
	\begin{abstract}
		In this paper, we study compact linear relations and establish that the spectrum of compact linear relations consists of only point spectrum and the zero set under suitable conditions.  
	\end{abstract}

	\maketitle
	\section{Introduction}
	Linear relations are generalization to linear operators.  The main difference is that a linear operator is a map between two linear spaces but is single-valued while a linear relation is a map between two linear spaces that can be multi-valued. In 1932, J. von Neumann introduced linear relations in connection to the non-densely defined linear differential operators, see \cite{von}. Later in 1961, Arens  developed some theory of linear relations in \cite{Arens}. Since then, there have been tremendous studies on the theory of linear relations \cite{codd,CD, DD,HSS,HSD, HSDS} and has been very useful tools in many applications \cite{Remling}. These linear relations are also  termed as linear subspaces and multivalued linear operators. There are several papers in the literature, see the monograph \cite{Cross} as an example and citation there in, which deal with the theory of multivalued linear operators. The generalization of the theories from linear operators to linear relations have been an interest of several researchers. For example, in \cite{RSX}, the authors study the boundedness and closedness of linear relations and show the Closed Graph Theorem for linear relations in Banach spaces and in \cite{AS},  the author have proved von Newmann's theorem for linear relations on Hilbert spaces.
	
	The goal of this paper is to discuss compact linear relations, obtain some of their properties, and generalize the spectral theorem for compact linear relation. Some of the studies on spectral theory of linear relations can be found in  \cite{KA, BK, BZ, GV}. We extend the theory for compact linear relations. In particular, we show that the spectrum of compact linear relations consists only of the point spectrum and the zero set under suitable conditions. \\

	Let $X$, $Y$ and $Z$ be linear spaces over the set of complex numbers $\C$.  A linear relation $T$ from $X$  to $Y$ is a multivalued map such that its graph
	\begin{equation*}
	G(T) = \{(x,y)\in X\times Y \  | \  y\in T(x)   \}
	\end{equation*} 
is a subspace	of the product space $X \times Y$. The set $T(x)= \{y\in Y : (x,y) \in G(T)\}\subseteq Y$ is the image of $x\in X$ under $T$. It is clear that the set $T(0)$ is multivalued part of $T$. In other words, for any  $x\in X,\ \  T(x) $ can be written in the form $$ T(x)= \{y\} +T(0). $$\\

The set of all linear relations from $X$ to $Y$ is denoted by $LR(X, Y)$. We  write  $LR(X)$ for $LR(X, X)$. The identity relation $I$ on $X$ is $
		I = \{(x,x): x \in X  \}
	$ which is obviously an identity operator on $X$.  The domain and the range of $T$ are denoted respectively as 
	\begin{equation*}
	D(T) = \{x \in X :  T(x) \neq \emptyset \}\textrm{ and } R(T) = \{y \in Y:  y\in T(x) \text{  for some  } x\in D(T) \}\,.
	\end{equation*}
 The kernel  of $T$ is a subspace of $D(T)$ 
	\begin{equation*}
	N(T) = \{x \in D(T) :  T(x) = T(0) \}\,, 
	\end{equation*}
	and the cokernel of $T$ is a quotient space $ Y/R(T).$

	 The inverse of $T$ is a relation denoted by $ T^{-1}$ such that
	\begin{equation*}
G(T^{-1}) = \{(y,x) : (x,y) \in G(T) \}\,.
	\end{equation*}
	
	The image of $M\subseteq X$ under $T$ is defined to be the set 
	\begin{equation*}
	T(M) = \bigcup_{x\in M} T(x)\,.
	\end{equation*} 
	A relation $T$ is injective if $N(T)=\{0\}$ and surjective if for every $y\in Y$ there is a $x\in D(T) $ such that $y \in T(x)$. Note that a relation $T$ is an operator if  $T(0)= \{0\}$ and $T^{-1}$ is an operator if $T$ is injective. A relation $T$ is called closed relation if its graph $G(T)$ a closed subspace of $ X\times Y$. The completion $\widetilde{T}$ of $T$ is defined by $G(\widetilde{T}):= \widetilde{G(T)}\subseteq \widetilde{X}\times \widetilde{Y},$ where $\widetilde{X}$ and $ \widetilde{Y}$ are completion of $X$ and $Y$ respectively.
	
	Suppose $T,T_1,T_2\in LR(X,Y)$ and $S\in LR(Y,Z)$ and $a$ a scalar, we define $ aT, \ T_1+ T_2, \ ST $  by the graphs
	\begin{equation*}
	\begin{aligned}
G(aT) & = \{(x,ay): (x,y) \in G(T) \}\\ 
	G(T_1+ T_2 )& =\{(x, y_1 +y_2): (x,y_1) \in G(T_1),\,\, (x,y_2) \in G(T_2) \} \\ 
	G(ST) & = \{(x, z) : (x, y) \in T, \,\, (y, z) \in S \}\,.
	\end{aligned}
	\end{equation*}
	
	Let $X, A$ be two linear spaces over $\C$.  We call $(X, A)$ a dual pair if there exists a mapping $\langle\cdot,\cdot\rangle:X\times A\to\C $ such that 
	\begin{equation*}
	\langle x,\lambda a\rangle = \bar{\lambda}\langle x, a\rangle
	\end{equation*}
	for all $x\in X$, $a\in A$ and scalars $\lambda$.  Let $(X, A)$ and $(Y, B)$ be two dual pairs over $\C$ and $T \in LR(X,Y) $ then $(X \times Y, A\times B)$ is a dual pair with
	\begin{equation*}
	\langle (x,y),(a,b)\rangle = \langle x,a\rangle +\langle y,b\rangle\,.
	\end{equation*}
	The adjoint of $T$ is a relation $T^*\in LR(B,A)$ given by 
	\begin{equation*}
	G(T^*) = \{(b,a): \langle x,a\rangle = \langle y,b\rangle, \,\textrm{ for all } (x,y)\in G(T) \}\,.
	\end{equation*}
	
	The following proposition is combination of propositions, lemmas and theorems due to Arens.  The proofs can be found in \cite{Arens}.
	
	\begin{prop}
		Let $T:X\to Y$ be a linear relation.  (1) $T\subseteq T^{**}$ and $T = T^{**}$ if and only if T is closed.  (2) If $S\subseteq T$, then $T^* \subseteq S^* $.  (3) $(\lambda T)^* = \bar{\lambda}T^*$. (4) $(T^{-1})^* = (T^*)^{-1}$ provided $T^{-1}$ exists. (5) $N(T^*) = R(T)^{\perp}$. (6) $T^*$ is single valued if and only if $D(T)$ is dense. (7) $T^*$ is closed. (8) If $S:X\to Y$ is also a linear relation then $(S+T)^* \supseteq S^*+ T^*$ and if $D(S^*)=B$ and $D(S) \supseteq D(T)$ then $(S+T)^*= S^* + T^*$.  (9) If $T\in LR(X,Y)$ and $S\in LR(Y,Z)$ and $(Z,C)$ is a dual pair then $(ST)^* = T^*S^*$ provided that $D(S^*) = C$ or $R(T^*) = A$.
	\end{prop}

  To obtain an operator from $T$ consider the quotient space $  Y/\overline{T(0)}.$ If $Y$ is a normed space, then so is the quotient space $Y/\overline{T(0)}$ with the norm defined by
	\begin{equation*}
	\|[y]\| = \textrm{inf}\{\|y - u\| : u \in \overline {T(0)} \} \,.
	\end{equation*}
	
	The quotient map $Q_T : Y\to Y/\overline{T(0)}$ is defined by $Q_T(y) = [y].$ Define the map $T_s$,
	\begin{equation*}
	T_s : D(T)\to Y/\overline{T(0)}
	\end{equation*}
	by $T_s = Q_T T$ so that $T_s(x) = [y]$ for all $(x,y)\in G(T) $. Then $T_s$ , called  an operator part of $T$, is an operator on $D(T)$  to $Y/\overline{T(0)}$. The norm of $T(x)$ for $x\in D(T)$ and the norm of $T$ are defined by respectively $\| T(x)\|= \|T_s(x)\|$ and  $\|T\| = \| T_s\|$. Note that this is a semi-norm since $\|T\| = 0$ does not imply that $T = 0.$  A relation $T\in LR(X,Y)$ is bounded if $\|T\| \leq M$ for some $M>0 .$
	The set of all bounded linear relations from $X$ to $Y$ is denoted by $BR(X,Y).$ We write $BR(X)$ for $BR(X, X).$
		
	It is shown in \cite{Cross} that $\|T^*\| \leq \|T\|$ for any $T \in LR(X, Y)$. Therefore we have the following proposition.

\begin{prop}\label{T_bnded_T_star_bnded}
	If $T\in BR(X,Y)$, then $T^* \in BR(Y^*, X^*)$ where $(X\times Y, X^*\times Y^*) $ is a dual pair.
\end{prop}

	\section {Spectrum of closed linear relations}
	
	 In this section, we discuss some spectral theory of linear relations and establish some relationship with so called the Arens' spectrum. We adopt the definition of resolvent set of a linear relation on a normed space by R. Cross from \cite{Cross}.	 Let $X$ be a complex normed space and  $T \in LR (X)$,  $\lambda \in \C$ write $T_{\lambda} = (\lambda I - \widetilde{T})^{-1}.$
	 
	\begin{defin} The {\it resolvent set}, $\rho(T)\subseteq\C$, of $T$  is defined as the set
		\begin{equation}\label{as}
			\rho(T) = \{\lambda \in \C :T_{\lambda} \text{  everywhre defined and single valued}  \}.
		\end{equation}
	
	\end{defin}
	
	Note that if $T$ is a closed relation on a Banach space $X$, the resolvent set is $$ \rho (T) = \{\lambda 
\in \C : (\lambda I - T)^{-1} \in B(X)\}$$

Thee {\it spectrum}, $\sigma(T)$, of $T$ is defined to be $\sigma(T)= \C\setminus\rho(T)$.

A complex number $\lambda\in\C$ is called an {\it eigenvalue} of   $T$ if $N(\lambda I -T)\neq \{ 0 \}.$ That is, there exists a nonzero $x\in X$ such that $ \lambda x \in T(x)$.   The set of all eigenvalues of $T$ is called the {\it point spectrum}, denoted by $\sigma_p(T)$.  It follows from the definition of the spectrum that $\sigma_p(T)\subseteq\sigma(T)$.\\

From now on, we consider $X$ a complex Banach space.  We can also define the {\it continuous spectrum}, $\sigma_c(T)$ and the {\it residual spectrum}, $\sigma_r(T)$, similar to the classical definition.  In fact we have:
	\begin{enumerate}
	\item $\sigma_p(T) = \{\lambda\in\C : N( \lambda I-T) \neq\{0\}\}$.
	\item $\sigma_c(T) = \{\lambda\in\C: N( \lambda I-T)=\{0\},\,\overline{R( \lambda I-T)} = X\}$.
  \item $\sigma_r(T) = \{\lambda\in\C: N( \lambda I-T)=\{0\},\,\overline{R(\lambda I-T)} \neq X\}$.
	\end{enumerate}
	It also follows that $\sigma_r(T)\subseteq\sigma(T)$ and $\sigma_c(T)\subseteq\sigma(T)$.  Moreover, the point, residual and continuous spectrum are pairwise disjoint and the union of all three form the spectrum.  Consequently, if $\lambda\in\rho(T)$, then $( \lambda I-T)^{-1}$ is a bijection and hence its inverse, $ \lambda I -T $ is a bijection on its range with the property that
\begin{equation*}
( \lambda I -T)( \lambda I -T )^{-1} = I_X \quad\textrm{and}\quad ( \lambda I-T)^{-1}( \lambda I-T) = I_{D(T)}\,,
\end{equation*}
where $I_X$ and $I_{D(T)}$ are the identity relations on $X$ and $D(T)$ respectively.  So, the above doesn't imply that $T-\lambda I$ is necessarily bounded.

	\begin{theo}
		Assume that  $T$ is a closed and bijective linear relation on a Banach space $X$ with $\|\lambda T^{-1}\|< 1$ for $|\lambda |>1$. Then $\lambda\in\rho(T)$ and 
		\begin{equation*}
			( \lambda I-T)^{-1} =  -\sum_{n=0}^{\infty} \lambda^n(T^{-1})^{n+1}\,,
		\end{equation*}
		where the series converges in operator norm.
	\end{theo}
	
	\begin{proof}
		Claim: $( \lambda I -T)$ is a bijection.  Suppose $ 0\in( \lambda -T)(x)$.  Then $\lambda x \in T(x)$ and hence $  x \in T^{-1}(\lambda x)$. 
		We have the following estimate:
		\begin{equation*}
			\| T^{-1}(\lambda x)\| \leq \| T^{-1}\|\|x\| < \frac{1}{|\lambda| }\|x\|\,.
		\end{equation*}
		Hence
		\begin{equation*}
			\|x\| < \frac{1}{|\lambda | }\|x\|\,.
		\end{equation*}
		Since $|\lambda|> 0,$ we must have $x = 0$. Therefore, $( \lambda I-T)$ is injective. Next we show that $R( \lambda I-T) = X$.  Indeed, suppose there is a $v\in X$, we must show that there is $ y \in T(x)$ for some $x\in X $ such that $v =  \lambda x -y$.  To ensure the existence of such a $  y $, define a map $L:X\to X$ by 
		\begin{equation*}
			L(y) =  -v + \lambda x  =  -v + \lambda T^{-1}y
		\end{equation*}
		for every $ y \in T(x) $.  Notice $L$ is well-defined and linear.  We have the following estimate:
		\begin{equation*}
			\|L(y) - L(\tilde{y})\| = \|\lambda T^{-1} y -\lambda T^{-1} \tilde{ y } \|\leq \|\lambda T^{-1} \| \| y - \tilde{y}\|\,.
		\end{equation*}
		Since $\|\lambda T^{-1}\|< 1$, $L$ is a contraction mapping on $X$. Thus by Banach contraction mapping theorem, there exists unique $y$ such that $L(y) = y$, that is $y = -v +\lambda T^{-1}y$. So $v  =   \lambda x-y$.  Next we compute the von Neumann series for $( \lambda-T)^{-1}$. Factoring out $T$ on the right of $ \lambda I -T$ we have 
		\begin{equation*}
		 \lambda I -T = ( \lambda T^{-1} -I )T\,,
		\end{equation*}
		and thus 
		\begin{equation*}
			( \lambda I -T)^{-1} = T^{-1}( \lambda T^{-1} -I )^{-1} = -  T^{-1}( I- \lambda T^{-1}  )^{-1}\,.
		\end{equation*} 
		Since $\lambda T^{-1}$ is a bounded linear operator, $( I - \lambda T^{-1})^{-1}$ has a von Neumann series, namely
		\begin{equation*}
			(I - \lambda T^{-1})^{-1} = \sum_{n=0}^{\infty} (\lambda T^{-1})^n =\sum_{n=0}^{\infty} \lambda^n ( T^{-1})^n\,.
		\end{equation*}
		Therefore,
		\begin{equation*}
			( \lambda I -T)^{-1} = -\sum_{n=0}^{\infty} \lambda^n ( T^{-1})^{n+1}\,.
		\end{equation*}
		It follows that $( \lambda I -T )^{-1}$ is a bounded linear operator and hence $\lambda\in\rho(T)$ completing the proof.
	\end{proof}

	The spectrum of a linear relation $T \in LR(X)$ can also be described through spectrum of linear operator $J : X \rightarrow X/X_0$ where $X_0$ is a closed subspace of $X$, see \cite{GV}. \\ 
	
	Let $X$ be a Banach space and $X_0$ is a closed subspace of $X$. Suppose 
	\begin{equation*}
		J: D(J) \subseteq X \rightarrow X/X_0
	\end{equation*}
	is a linear operator, then the Arens's resolvent set of $J$ is defined by
	
	\begin{equation*}
		\rho_A(J)= \{ \lambda \in \C : (\lambda Q -J)^{-1} \in B(X/X_0 , X) \}\, .
	\end{equation*}
	
	Here $Q$ is a natural quotient map from $X$ onto $X/X_0$. The Arens's spectrum is defined by 
	\begin{equation*}
		\sigma_A(J)= \C \setminus \rho_A(J)\,.
	\end{equation*}
	 If $T \in LR (X)$ is closed,  then clearly $T(0)$ is closed and the  the operator $T_s$ is also closed. For $X_0 = T(0)$,  $J=T_s$ and  $I_s = Q_T I $  the operator $\lambda I_s - T_s$ is closed   for any $\lambda \in \C .$  Indeed, the Arens's resolvent set coincides with the resolvent set defined in \eqref{as}.
	\begin{prop}\label{prop1} Let $X$ be a Banach space and $T \in BR(X)$ be closed. Then,
		\begin{equation*}
			\sigma_A(T_s) = \sigma (T)\,.
		\end{equation*}
	\end{prop}
	
	\begin{proof} 
		To show $\sigma_A(T_s)=\sigma(T)$, we show their respective resolvent sets are equal.	 First let $\lambda \in \rho(T),$ then $ (\lambda I - T)^{-1} \in B(X)$.   We show that $\lambda I_s- T_s$ is bijective. Let $x_1, x_2 \in D(T)$ such that $x_1 \neq x_2$.  Suppose $ (\lambda I_s - T_s)(x_1)= (\lambda I_s - T_s)(x_2)$. Then \begin{align*}  & \lambda [x_1] -[y_1]=  \lambda[x_2] -[y_2]  \end{align*}
		
		for some $ y_1 \in T(x_1), \  y_2 \in T(x_2) $  and so $ [\lambda(x_1-x_2)] = [y_1-y_2] .$ It follows that \begin{equation}\label{eq 1} \lambda(x_1-x_2) - (y_1-y_2) \in T(0) .\end{equation} On the other hand \begin{equation} \label{eq 2} (  y_1-y_2) \in T(x_1-x_2) . \end{equation}
		
		Adding equations \eqref{eq 1} and \eqref{eq 2} yields
		\begin{equation*}
		 \lambda(x_1-x_2) \in T(x_1-x_2) \,.
		\end{equation*}
		Thus $(x_1-x_2) \in N (\lambda I -T)$.  However $N(\lambda I -T)= \{0\}$ since $\lambda I -T$ is injective.  Hence $x_1 = x_2$ which is a contradiction. Therefore 
		\begin{equation*}
			(\lambda I_s - T_s)(x_1) \neq (\lambda I_s - T_s)(x_2)\,.
		\end{equation*}
		
		Next we show $\lambda I_s -T_s$ is surjective.  Since $\lambda \in \rho(T)$ and $T$ is a closed relation,  $R( \lambda I -T) = X$, it follows that $R(\lambda I_s- T_s)= X/T(0)$ and hence surjective.  Since $T$ is closed and bounded, $\lambda I -T$ is closed and bounded and hence $(\lambda I -T)^{-1}$ is bounded, it follows that $\lambda I_s- T_s$ is closed and  $(\lambda I_s- T_s)^{-1}$ is bounded. Thus we have 
		
		\begin{equation} \label {eq 3}
			\rho(T) \subseteq \rho_A(T_s) \,.
		\end{equation}
		Secondly, let $\lambda\in\rho_A(T_s)$. Then $(\lambda I_s- T_s)^{-1} \in B(X/T(0), X)$. We show that $(\lambda I -T)^{-1} \in B(X)$. Let $u \in X$. Then,  
		\begin{equation*}
			x = ( (\lambda I_s- T_s)^{-1}) [u] \in X \textrm{ and }(\lambda I_s- T_s)(x) = [u]\,.
		\end{equation*}  
		This implies 
		\begin{equation*}
			\lambda [x]- [y]= [u]
		\end{equation*}
		for some $ y \in T (x) $. It follows that there is some $y \in T(x)$ such that $\lambda x -y =u .$ That is $u \in R(\lambda I -T)$ and so $R( \lambda I -T) = X$, hence $\lambda I-T$ is onto.  
		
		Next suppose $x\in N(\lambda I -T)$, then we must have $\lambda x \in T(x) $. This implies $T_s(x) = \lambda [x] $ so that $x \in N(\lambda I_s- T_s ) = \{0\}$. Thus $x=0$ which implies that $\lambda I -T$ is injective.  Therefore $\lambda I -T$ is bijective and so $(\lambda I -T)^{-1}$ exists as an operator. 
		
		Notice that $(\lambda I -T)^{-1}$ is bounded because $(\lambda I_s- T_s )^{-1}$ is bounded. Hence $\lambda \in \rho (T)$. Thus,
		
		\begin{equation} \label {eq 4}
			\rho_A(T_s) \subseteq \rho(T)
		\end{equation}
		
		Thus the proof is complete. 
	\end{proof}
	
	Next, we define 
	\begin{equation}\label{t_tilde_def}
		 \mathbb{T} : X/T(0) \rightarrow X/T(0) \textrm{ by }  \mathbb{T} [x] = y+T(0) =[y] \textrm{ for some } y \in T(x)\, .
	\end{equation}
	
	Notice that $ \mathbb{T}$ is a single valued linear operator if $T(T(0)) \subseteq T(0)$ and $ D(T) = X$. For if $[x]= [u]$, then $ \mathbb{T}([x])= y+T(0)$ for some $y \in T(x)$ and $ \mathbb{T}([x])= y'+T(0)$ for some $y' \in T(u)$.	Since $x-u \in T(0)$ and $T(T(0))\subseteq T(0)$, we have $T(x-u) \subseteq T(0)$.  Since $T$ is linear $Tx-Tu \subseteq T(0)$. This implies $ y-y' \in T(0)$ and $y +T(0) = y'+ T(0)$.  Since $ \mathbb{T}$ is linear, it follows that $ \mathbb{T}$ is a single valued linear operator. Then for any $x\in X$ we have, 
	\begin{equation*}
		 \mathbb{T}([x])=T_s(x)\,.
	\end{equation*}
	
	\begin{prop}\label{prop2}
		Let $T\in BR(X)$ be closed,  $T(T(0)) \subseteq T(0)$ and $ D(T) = X$, and $ \mathbb{T}$ be defined as in equation (\ref{t_tilde_def}), then 
		\begin{equation*}
			\sigma_p(T) =  \sigma_p( \mathbb{T})\,.
		\end{equation*}
	\end{prop}
	
	\begin{proof}
		Let $\lambda \in \sigma_p(T)$. There exists $ x\in X$ with $x\neq0$ such that $ \lambda x \in T(x)$. By definition of $ \mathbb{T}$ we have $ \mathbb{T}([x]) = \lambda [x]$. So $\lambda \in \sigma_p( \mathbb{T})$. Similarly, let $\lambda \in \sigma_p( \mathbb{T})$. Then $ \mathbb{T}([x]) = \lambda [x]$. Suppose $ y  \in T(x) $, then $ \mathbb{T}([x]) =  [y]$ so that $ [y] =  \lambda [x] $ and therefore $\lambda x -y \in T(0)$. This says that $\lambda x -y \in T(0)$. But $ y \in T(x)$ implies that $\lambda x \in T(x) $. Thus $ \lambda \in \sigma_p(T)$.  This completes the proof.
	\end{proof}
	
	\begin{prop}\label{prop3} 
		Let $T\in BR(X)$ be closed and $ \mathbb{T}$ be defined as in equation (\ref{t_tilde_def}), then 
		\begin{equation*}
			\sigma_A (T_s) =  \sigma ( \mathbb{T})\,.
		\end{equation*}
	\end{prop}

	\begin{proof} We will show that $\rho (  \mathbb{T}) = \rho_A (T_s)$.
		
		Let $\lambda\in\rho( \mathbb{T})$, then $(\lambda I - \mathbb{T})^{-1} \in B (X/T(0))$. Since $R(\lambda I_s -T_s) = R(\lambda I - \mathbb{T})= X/T(0)$ we have that $( \lambda I_s -T_s)$ is surjective. Suppose $x_1, x_2 \in X$ such that $x_1 \neq x_2$. If 
		\begin{equation*}
			(\lambda I_s -T_s)(x_1) = (\lambda I_s -T_s)(x_2)\,,
		\end{equation*}
		then $\lambda [x_1] - [y_1] = \lambda [x_2] - [y_2]$ for some $ y_1\in T(x_1)$ and $ y_2\in T(x_2)$. It follow that $\lambda [x_1 -x_2] = [y_1 -y_2]$. This means that
		\begin{equation*}
			\lambda(x_1-x_2) -(y_1-y_2) \in T(0)
		\end{equation*}
	  Given $y_1 \in T(x_1)$ and $y_2 \in T(x_2)$, linearity of $T$ implies that $ ( y_1 -y_2) \in T(x_1 -x_2)$. Therefore 
		\begin{equation*}
		 \lambda(x_1-x_2) ) \in T	(x_1 -x_2) \,.
		\end{equation*}
		This shows that $\lambda \in \sigma_p(T) =  \sigma_p( \mathbb{T})$ by Proposition \ref{prop2}.  However this is a contradiction $\lambda\in\rho( \mathbb{T})$. Thus 
		\begin{equation*}
			(\lambda I_s -T_s)(x_1) \neq ( \lambda I_s -T_s)(x_2)
		\end{equation*}
		proving that $\lambda I_s -T_s$ is injective.  We have that $(\lambda I_s -T_s)^{-1}$ is a well-defined operator. The boundedness of $(\lambda I_s -T_s)^{-1}$ follows from the boundedness of $(\lambda I -   \mathbb{T})^{-1}$. 
		
		We show the other containment.  Let $\lambda\in\rho_A (T_s)$ then $(\lambda I_s -T_s)^{-1} \in B(X/T(0), X)$ and suppose that
		\begin{equation*}
			(\lambda I -  \mathbb{T})([x_1]) = (\lambda I -   \mathbb{T}) ([x_1])
		\end{equation*}
		for some classes $[x_1]$ and $[x_2]$ in $X/T(0)$.  Then 
		\begin{equation*}
			(\lambda I_s -T_s)(x_1) = ( \lambda I_s -T_s)(x_2)
		\end{equation*}
		for $x_1$ and $x_2$ in $X$. Since $\lambda I_s -T_s$ is injective, we have that $x_1 =x_2$. This implies that $[x_1] = [x_2]$.  It follows that $\lambda I -   \mathbb{T}$ is injective. Again, since 
		\begin{equation*}
			R(\lambda I_s -T_s) = R(\lambda I - \mathbb{T})= X/T(0)\,.
		\end{equation*}
		We have that $\lambda I -  \mathbb{T}$ is surjective. Therefore $\lambda I -   \mathbb{T}$ is bijective. The boundedness of  $(\lambda I -   \mathbb{T})^{-1}$ follows from the boundedness of $(\lambda I_s -T_s)^{-1}$. Hence $\lambda\in\rho( \mathbb{T})$. This proves $\rho( \mathbb{T}) = \rho_A (T_s)$ completing the proof.
	\end{proof}
	
	\section{ Compact linear relations and their spectral properties }
	\begin{defin}
		Let $X$ and $Y$ be normed spaces. A linear relation $T: X \to Y$ is called a {\it compact} if for every bounded subset $M$ of $X$, the image $T_s(M)$ is precompact, that is, $\overline{T_s(M)}$ is a compact subset of $ Y/\overline{T(0)}$.
		\end{defin}  
	This definition coincides with the definition of compact operators. 
	We denote the set of all compact linear relations by $KR(X,Y)$ and use $KR(X)$ when $X=Y$. Note that since compact subsets in a Banach space are necessarily bounded, any compact relation is bounded. Therefore, we have $ KR(X, Y) \subseteq BR (X, Y)$. The goal of the following propositions is to extend some of the classical functional analytic results about compact operators to compact linear relations.\\
	
  In this section we will discuss compact relations and their spectral properties.  In \cite{MH}, the authors describe so called the congruent spectrum for compact relations. Some of the results in the paper may be parallel but have different perspective.  
	
	\begin{prop} A linear relation $T\in LR(X,Y)$ is compact if and only if every sequence $\{x_n\}\in D(T)$, with $\|x_n\| <1$ has a subsequence $ \{x_{n_k}\} $ for which $T_s (x_{n_k})$ converges in $ Y/\overline{T(0)}.$
	\end{prop}
	\begin{proof}
	Clear from the definition.
	\end{proof}
	\begin{prop}
	Let $X$ and $Y$ be Banach spaces.  If $T\in KR(X,Y)$ then $T^* \in KR(Y^*,X^*)$. 
	\end{prop}
	
	\begin{proof}
		Let $T\in KR(X,Y)$, $B_X$ be the unit ball in $X$ and $B_{Y^*}$ be the unit ball in the dual space of $Y$, $Y^*$. Then 
		\begin{equation*}\| T^*f\| = \sup_{\substack{x\in B_X\\ x\in D(T)}} \|T_s^*f (x)\|=\sup_{\substack{x\in B_X\\ x\in D(T)}} \| f(T_s (x))\|= \sup_{y\in \overline{T(B_X)}} \|f(y)\|\,.
		\end{equation*}
		Define $A = T^*(B_{Y^*})$ and embed $A\subset X^*$, with $X^*$ the dual space to $X$, into $C(\overline{T(B_X)})$ via 
		\begin{equation*}
		u:A\to C(\overline{T(B_X)}),\,\, u(T^*f) = f|_{\overline{T(B_X)}}
		\end{equation*}
		for $f\in B_{Y^*}$.  Thus $\|T^*f\|_{X^*} = \|f|_{\overline{T(B_X)}} \|_{C(\overline{T(B_X)})}$, for all $f \in Y^*$. 
		
		By Proposition \ref{T_bnded_T_star_bnded}, if $T$ is bounded then $T^*$ is bounded and $\|T\|= \|T^*\|$ so that 
		\begin{equation*}
		\|T^*\| \leq \|T^*\|\|f\|_{Y^*} \leq \|T\|\,.
		\end{equation*}
		Therefore, $u(A)$ is uniformly bounded in $C(\overline{T(B_X)})$. Given $y_1, y_2 \in\overline{T(B_X)},\,f\in B_{Y^*}$, we have
		\begin{equation*}
		|f|_{\overline{T(B_X)}} (y_1)- f|_{\overline{T(B_X)}} (y_2) | \leq | f(y_1-y_2)| \leq \|f\|_{Y^*} \|y_1 -y_2\| \leq \|y_1 - y_2\|\,.
		\end{equation*}
		Thus, $u$ is equicontinuous and hence by Arzela and Ascoli $A$ is precompact.  This completes the proof. 
	\end{proof}
	Moreover, it is shown in \cite{Cross} that if $T$ is closed in above proposition then $T\in K(X,Y)$ if and only if $T^* \in K(Y^*, X^*). $  
	\begin{prop}\label{ideal}
		Let $T\in KR(X)$ and $S\in BR(X)$ be such that $TS$ and $ST$ are well-defined compositions, that is $R(T)\subseteq D(S)$ and $R(S)\subseteq D(T)$.  Then, $TS\in KR(X)$ and $ST\in KR(X)$.
	\end{prop}
	
	\begin{proof}
		Let $U$ be a bounded subset of $D(S)$. Since $S\in BR(X)$, we have $S(U)$ is a bounded set. Therefore $(TS)(U) = T(S(U))$ which is a precompact set since $T\in KR(X)$.   Now let $U$ be a bounded set in $D(T)$.  By definition, since $T(U)$ is precompact and $S$ is bounded and hence continuous, it follows that $(ST)(U)$ must be precompact.
	\end{proof} 
	This says that $KR(X)$ behaves like an ideal in $B(X)$. 
	
	\begin{prop}
		Let $X$ be a Banach space and $T\in KR(X)$ such that $ D(T)$ and $T(0)$ are closed. Then for any $\lambda\in\C $, $N(\lambda I-T)$ is a closed subspace of $X$.
	\end{prop}
	
	\begin{proof}
		Let $x\in\overline{N(\lambda I-T)}$. Then there exists $x_n\in N(\lambda I-T) $ such that $x_n$ converges to $x$. Since $x_n\in N(\lambda I-T)$, we have $(x_n, 0)\in\lambda I-T$, and thus $(x_n,\lambda x_n)\in T$.  Since $x_n$ converges to $x$, it follows that $(x_n, \lambda x_n)$ converges to $(x, \lambda x)$. By assumption $D(T)$ and $T(0)$ being closed, $T$ is closed, and hence $(x, \lambda x) \in T$. It follows that $(x,0)\in\lambda I - T$ and hence $x\in N(\lambda I - T)$. Therefore $N(\lambda I - T)$ is a closed subspace of $X$.
	\end{proof}

\begin{theo}
Let $X$ be a Banach space and  $T \in LR (X)$ be a closed relation with $D(T)= X$ and $T(T(0)) \subseteq T(0).$ Then $T$ is compact if and only if $ \mathbb{T}$ is compact.
\end{theo}

\begin{proof}
Suppose $T$ is compact. Let $\{[x_n]\}$ be a bounded sequence in $X/T(0)$. Then there exists $M>0$ such that $\|[x_n]\| \leq M$ for all $n$. Then by the approximation property of the infimum, for each $n$ there exists $u_n \in T(0)$ such that 
\begin{equation*}
\|x_n -u_n\| \leq \|[x_n]\| + \frac{1}{2^n} \leq M + \frac{1}{2^n}\,.
\end{equation*}
Hence $\{x_n-u_n\}$ is a bounded sequence in $X$.  Since $T$ is compact, there exists a subsequence $\{x_{n_k}-u_{n_k}\}$ such that $T_s(x_{n_k}-u_{n_k})$ converges. Since $T_s = QT$ and $T(u_{n_k})\in T(0)$, it follows that
\begin{equation*}
T_s(x_{n_k}-u_{n_k}) = T_s (x_{n_k})\,.
\end{equation*}
Moreover $T_s(x_{n_k})=  \mathbb{T}([x_{n_k}])$, hence there is a subsequence $\{[x_{n_k}]\}$ such that $ \mathbb{T}([x_{n_k}])$ converges. Therefore $ \mathbb{T}$ is compact. 

Conversely, suppose $ \mathbb{T}$ is compact. Let $\{x_n\}$ be a bounded sequence in $X$. Then $\{[x_n]\}$ is a bounded sequence in $X/T(0)$ since $0\in T(0)$ and 
\begin{equation*}
\|[x_n]\| \leq \|x_n-0\|= \|x_n\|\,.
\end{equation*}
Since $ \mathbb{T}$ is compact, there exists a subsequence $\{[x_{n_k}]\}$ such that $  \mathbb{T}([x_{n_k}])$ converges. For each $n_k$, 
\begin{equation*}
 \mathbb{T}([x_{n_k}]) =  T_s (x_{n_k})\,,
\end{equation*}
thus there is a subsequence $\{x_{n_k}\}$ such that $T_s (x_{n_k})$ converges. Therefore $T$ is compact.
\end{proof}

\begin{theo}[The Spectral Theorem]
Let $T$ be a compact linear relation on a Banach space $X$ with $\dim{X}= \infty $ such that $D(T)= X$ and $T(T(0)) \subseteq T(0)$. Then
\begin{equation*}
\sigma(T) = \sigma_p(T)\cup\{0\}\,.
\end{equation*}
In particular, $0$ is an accumulation point for $\sigma_p(T)$.  Moreover $0\in\sigma_p(T)$ if $\sigma_p(T)$ is closed, $0\in\sigma_r(T)$ if $T$ is not onto and its range is not dense, and $0\in\sigma_c(T)$ if $T$ is not onto but has dense range.
\end{theo}
	
\begin{proof}
First notice that $0\in\sigma(T)$.  Otherwise $T$ is invertible and its inverse is bounded. Therefore by Proposition \ref{ideal}, $I=T^{-1}T\in KR(X)$ which is impossible. Since $T$ is compact, $ \mathbb{T}$ is compact and therefore by the Spectral Theorem for compact operators, 
\begin{equation*}
\sigma(\ \mathbb{T})\setminus\{0\} = \sigma_p ( \mathbb{T})\,.
\end{equation*}
By Propositions \ref{prop1}, \ref{prop2} and \ref {prop3}, we have
\begin{equation*}
\sigma(T) = \sigma_p(T)\cup\{0\}\,.
\end{equation*}
Moreover since $\sigma_p(T) = \sigma_p(\mathbb{T})$, and $ \mathbb{T}$ is a compact operator, $0$ is an accumulation point of $\sigma_p( \mathbb{T})$.
\end{proof}

The notion of essential spectrum of a bounded linear relation on a Banach space is a natural question. This has been considered by many authors, see the book \cite{Cross} for introduction and \cite{Alvarez, Wilcox2} for current research.  The essential spectrum of a linear relation on Banach space is similarly defined as for an operator. Let $X$ and $Y$ be Banach spaces and $T \in BR(X, Y) $ be a closed bounded linear relation.  $T$ is called a {\it Fredholm relation} if $T$ has finite dimensional kernel and cokernel and closed range.  The {\it essential spectrum} of $T$ is defined to be
\begin{equation*}
\sigma_{\textrm{ess}}(T) = \{\lambda\in\C : T-\lambda I \textrm{ is not Fredholm}\}
\end{equation*}
   It is well known that for a compact operator, $T$, we have $\sigma_{\textrm{ess}}(T)=\{0\}$.  So it is   natural to ask what is the essential spectrum for a closed compact linear relation on a Banach space $X$.   It is suspected that the essential spectrum of any closed compact linear relation is contained in $\{0\}$.  This suspicion is based on the fact that if $T$ is a compact operator on a Banach space $X$, then $T+I$ has finite dimensional kernel and cokernel.  This is typically shown by showing that the unit ball of the kernel of $T+I$ must be compact and hence by Riesz's Theorem, must be finite dimensional.  So, it seems plausible that it could happen for a closed compact linear relation on a Banach space.


\begin{thebibliography}{99}
		
		\bibitem{KA}
	 K. R. 	Acharya, {\it Self-Adjoint Extension and Spectral Theory of a Linear Relation in a Hilbert Space}, ISRN Mathematical Analysis, Art. ID 471640, 1 - 5, 2014.
		
	\bibitem{Alvarez}
	T. \'{A}lvarez, {\it On Essential Spectra of Linear Relations and Quotient Indecomposable Normed Spaces}, Rock. Mount. Jour. Math., 42, 1 - 15, 2012.
	
		\bibitem{Arens}
		 R. Arens, {\it Operational calculus of linear relations}, Pac. Jour. Math., 11, 9 – 23, 1961.
\bibitem{MH} M. Ayedi and H. Baklouti, {\it Congruent spectrum for compact linear relations,} Publ. Math. Debrecen 87/1-2 (2015), 191-207 DOI: 10.5486/PMD.2015.7138		
		\bibitem{BK}
	 A. G. 	Baskakov,  I. A. Krishtal, {\it On completeness of spectral subspaces of linear relations and ordered pairs of linear operators}, Jour. Math. Anal. Appl., 407, 157 - 178, 2013.
		
		\bibitem{BZ}
	 A. G.	Baskakov,  A. S. Zagorskii, {\it Spectral Theory of Linear Relations on Real Banach Spaces}, Math Notes, 81, 1 - 15, 2007. 
		
		\bibitem{codd}
	 E. A. 	Coddington, {\it  Extension Theory of Formally Normal and Symmetric Subspaces}, Mem. Amer. Math. Soc., 134, 1973.
		
		\bibitem{CD}
	 E. A.	Coddington, H. S. V. de Snoo,  { \it Positive selfadjoint extensions of positive symmetric subspaces}, Math. Z., 159, 203 - 214, 1978.
		
		\bibitem{Cross}
		 R. Cross, {\it Multivalued Linear Operators}, Marcel Dekker, New York, 1998.
		
		\bibitem{DD}
		  A. Dijksma,  H. S. V. de Snoo, {\it Self-adjoint extensions of symmetric subspaces},  Pac. Jour. Math., 54, 71 – 100, 1974.
		
		\bibitem{GV}
		 D.Gheorghe, F. H. Vasilescu, {\it Spectral theory for linear relations via linear operators}, Pac. Jour. Math., 255, 349 - 372, 2012.
		 
		 \bibitem{von}
		J. von Neumann, {\it \"Uber adjungierte Funktional-operatoren}, Ann. Math.,33(1932, 294-310) - 372, 2012.
		
		\bibitem{HSS}
		 S. Hassi,  A. Sandovici,  H. S. V. de Snoo, {\it Form sums of nonnegative selfadjoint operators}, Acta. Math. Hungar., 111, 81 - 105, 2006.
		
		\bibitem{HSD}
		 S. Hassi,  A. Sandovici,  H. S. V de Snoo, et. al., {\it  Extremal extensions for the sum of nonnegative selfadjoint relations}, Proc. Amer. Math. Soc., 135, 3193 - 3204, 2007.
		
		\bibitem{HSDS}
	 S. Hassi,  H. S. V.de Snoo, and  F. H. Szafraniec, {\it  Componentwise and canonical decompositions of linear relations}, Diss. Math., 465, 1 - 59, 2009.
		
		\bibitem{Remling}
		 C. Remling, {\it Spectral Theory of Canonical Systems}, De Gruyter Studies in Mathematics Series. Walter De Gruyter
		GmbH, 2018.
		
			
		\bibitem{RSX}
		 G. Ren,  Y. Shi, and  G. Xu, {\it Boundedness and Closedness of Linear Relations}, Lin. Multi. Alg., 66, 309 - 333, 2018.
		
		\bibitem{AS}
		 A. Sandovici, {\it Von Neumann’s theorem for linear relations}, Lin.and Multi. Alg., 2017.
		
		\bibitem{Wilcox}
	 D. L. 	Wilcox, {\it Multivalued semi-Fredholm operators in normed linear spaces,} Ph.D. thesis, University of Cape Town, 2002.
	
	\bibitem{Wilcox2}
	 D. L. Wilcox, {\it Essential Spectra of Linear Relations}, Lin. Alg. Appl., 462, 110 - 125, 2014.
	\end{thebibliography}
\end{document}